\documentclass{article}
\newtheorem{theorem}{Theorem}

\usepackage{psfig}
\def\ds{\displaystyle}
\def\aea{ & = & }
\def\ra{\rightarrow}

\def\f{\frac}

\newcommand\h[1]{\ensuremath{#1}}
\newcommand\beas[1]{\begin{eqnarray*} #1 \end{eqnarray*}}
\newcommand\bea[1]{\begin{eqnarray} #1 \end{eqnarray}}
\newcommand\lrs[1]{\left[ #1 \right]}
\newcommand\lrr[1]{\left( #1 \right)}
\newcommand\lrc[1]{\left\{ #1 \right\}}

\newcommand\g[1]{\Gamma\lrr{#1}}
\title{A Property of the Gamma Function at its Singularities}
\author{Anirudh Prabhu\\
West Lafayette Jr/Sr High School\\
West Lafayette, IN 47906\\
aprabhu@purdue.edu}
\date{}
\begin{document}
\maketitle
\begin{abstract}
The singularities of the \h{\Gamma} function, a meromorphic function on the complex plane,  are known to occur at the
nonpositive integers. We show, using Euler and Gauss identities, that for all positive integers \h{n}
and \h{k},
\beas{
\ds\lim_{z\ra 0} \ds\f{\g{nz}}{\g{z}} = \f 1 n; \hspace{0.4in}
\ds\lim_{z\ra -k} \ds\f{\g{nz}}{\g{z}} = \ds\f{(-1)^{k}\  \Gamma(k)}{n^2\  \Gamma(nk)}.} The above relations add to
the list of the known fundamental Gamma function identities.
\end{abstract}

\vspace{\baselineskip}

The Gamma function, which extends the factorial function  to the complex
plane, can be defined, following Euler and Weierstrass, as \cite{arfk}
\beas{
\Gamma(z) := \f1z \ds\prod_{n=1}^\infty \ds\f{\lrr{1+\f1n}^z}{1+\f z n} \hspace{0.2in}
}
The function is known to be meromorphic and its only singularities occur at the nonpositive integers
\cite{arfk}.  The \h{\Gamma} function can also be defined
in integral form for $Re(z)>0$ as
\beas{
\Gamma(z) = \int_0^\infty e^{-t} t^{z-1} dt
}
which can be extended to the region where \h{Re(z) \leq 0}, except for nonpositive integers, by
analytic continuation \cite{arfk}.  For other definitions of the \h{\Gamma} function
see \cite{arfk}.

The behavior of meromorphic functions near their singularities is a topic of
considerable interest in complex analysis. In the following we establish a property of
the \h{\Gamma} function at its singularities. First, we will establish the result stated in the
abstract for the singularity at \h{z=0}, and then extend the result to all the nonzero
singularities of the \h{\Gamma} function.
\begin{theorem}
For every positive integer \h{n},
\beas{\lim_{z\ra 0} \ds\f{\g{nz}}{\g{z}} = \f1n}
\end{theorem}
{\bf Proof:}  In the   Euler reflection formula for the $\Gamma$ function \cite{havi}
\bea{
\hspace{1.75in} \g \xi \g{1-\xi} = \ds\f\pi{\sin(\pi \xi)}\label{erf}
}
if we set \h{\xi=\ds\f k n}, for \h{1\leq k \leq n}, we obtain
\beas{
\sin\lrr{\ds\f{k\pi} n} \aea \ds \f{\pi}{\g{\ds \f k n}\g{\ds\f{n-k}n}}
}
which implies
\bea{
\ds\prod_{k=1}^{\left\lfloor n/2 \right\rfloor} \sin\lrr{\ds\f{k\pi} n} \aea \ds\prod_{k=1}^{\left\lfloor n/2 \right\rfloor}
\lrs{\f{\pi}{\g{\ds \f k n}\g{\ds\f{n-k}n}}}. \label{mainp}
}
The product on the left hand side of (\ref{mainp}) is related to the product of distances between
all pairs of \h{n} points that are uniformly distributed on a unit circle \cite{absa,aime}.  The following identity has been established using the Vandermonde determinant in \cite{absa} and more directly in \cite{plan} for positive integer \h{n}
\bea{\ds\prod_{k=1}^{n-1} \sin\lrr{\f{k\pi}n}= \f{n}{2^{n-1}}. \label{eq2}}
An outline of the proof of (\ref{eq2}), as reported in \cite{plan}, is presented in the Appendix.

Since \h{\sin\lrr{\f{k\pi}n} = \sin\lrr{\f{(n-k)\pi}n}} for $1\leq k\leq n$, and for even \h{n},
\h{\sin\lrr{\f{k\pi}n}=1}, for $k=\ds\f n 2$, we have
\bea{
\prod_{k=1}^{n-1} \sin\lrr{\f{k\pi}n} =
\lrs{\ds\prod_{k=1}^{\left\lfloor n/2 \right\rfloor} \sin\lrr{\ds\f{k\pi} n}}^2 \Rightarrow
\lrs{\ds\prod_{k=1}^{\left\lfloor n/2 \right\rfloor} \sin\lrr{\ds\f{k\pi} n}} = \lrs{\ds\f n{2^{n-1}}}^{\f12}
\label{sinf}
}
The denominator on the right hand side of (\ref{mainp}) can be rewritten as
\bea{
\prod_{k=1}^{\left\lfloor n/2 \right\rfloor} \g{\ds \f k n}\g{\ds\f{n-k}n} =
\left\{ \begin{array}{ll}
\ds\prod_{k=1}^{n-1} \g{\ds \f k n} & \mbox{\h{n} odd}\\
\\
(\ds \sqrt{\pi} )\ds\prod_{k=1}^{n-1} \g{\ds \f k n} & \mbox{\h{n} even}
\end{array}\right. \label{gam}
}
For even \h{n} there are two \h{\g{\f12}} factors on the left hand side above.
One of the them is included in the product while the other appears as \h{\g{\f12}=\sqrt\pi}.
Observing that for odd \h{n},  \h{\ds \pi ^{\lfloor n/2\rfloor}= \pi ^{n/2}\cdot \pi^{-1/2}}
and using (\ref{gam}) we have
\bea{
\ds\prod_{k=1}^{\left\lfloor n/2 \right\rfloor}
\lrs{\f{\pi}{\g{\ds \f k n}\g{\ds\f{n-k}n}}} \aea \ds\f{\pi^{(n-1)/2}}{\ds\prod_{k=1}^{n-1} \g{\ds \f k n}}
\label{rhs}
}
Next, consider the Gauss multiplication formula for the \h{\Gamma} function \cite{abst,emot}
\bea{
\hspace{1.75in} \ds\prod_{k=0}^{n-1} \g{z+\ds \f k n} = \lrr{2\pi}^{ \f{n-1}2}\ n^{ \f12 - nz}\  \g{nz} \label{g1}
}
which can be rewritten as
\bea{
\ds\prod_{k=1}^{n-1} \g{z+\ds \f k n} = \lrr{2\pi}^{ \f{n-1}2}\ n^{ \f12 - nz}\  \lrc{\ds\f{\g{nz}}{\g{z}}} \label{g2}
}
Taking the limit as \h{z\ra 0} on both sides we get
\bea{
\lim_{z\ra 0} \ds\prod_{k=1}^{n-1} \g{z+\ds \f k n} = \ds\prod_{k=1}^{n-1} \g{\ds \f k n} = \lrc{\lrr{2\pi}^{ \f{n-1}2}\ n^{ \f12}}\lrs{\lim_{z\ra 0}
\lrc{\ds\f{\g{nz}}{\g{z}}}} \label{eqw}
}
Using (\ref{mainp}), (\ref{sinf}), (\ref{rhs}) and (\ref{eqw}) we have
\beas{\lrs{\ds\f n{2^{n-1}}}^{\f12} = \ds\f{\pi^{(n-1)/2}}{\ds
\lrr{2\pi}^{ \f{n-1}2}\ n^{ \f12}\lrs{  \lim_{z\ra 0}
\ds\lrc{\ds \f{\g{nz}}{\g{z}}}}}}
or
\beas{
\lim_{z\ra 0}\ \
\ds \ds \f{\g{nz}}{\g{z}} =
\ds\f 1 n
}
as claimed \hfill \rule{7pt}{7pt}

\vspace{0.2in}

In the following we establish the result stated in the abstract
for all the nonzero singularities of the \h{\Gamma}
function.

\begin{theorem}
For all positive integers $n$ and $k$,
\beas{
\lim_{z\ra -k} \ds\f{\Gamma(n z)}{\Gamma(z)} = \ds\f{(-1)^k \ \Gamma(k)}{n^2 \ \Gamma(nk)}
}
\end{theorem}
{\bf Proof:}  In the Euler reflection formula (\ref{erf}) we set \h{\xi=-w+\f r n}, for \h{1\leq r\leq n-1} to obtain
\beas{
\Gamma\lrr{-w+\f r n} \Gamma\lrr{1+w-\f r n} = \ds\f{\pi}{\sin\lrr{\f{\pi r}n - \pi w}}
}
We observe that for a positive integer \h{k},
\bea{
\lim_{w\ra k} \lrs{\ds\f{\pi}{\sin\lrr{\f{\pi r}n - \pi w}}} = \ds\f{(-1)^k\ \pi}{\sin\lrr{\f{\pi r}n}}
\label{eq10}
}
Using the above equation and (\ref{eq2}) we get, for positive integer \h{k},
\bea{
\lim_{w\ra k} \prod_{r=1}^{n-1} \Gamma\lrr{-w+\f r n}\Gamma\lrr{1+w-\f r n} \aea \ds\f{(2\pi)^{n-1} (-1)^k}{n} \label{b1}
}
Using Gauss' multiplication identity (\ref{g1}) and rewriting it as in (\ref{g2}), except using \h{-w} instead of
\h{z}, we get
\bea{
\ds\prod_{r=1}^{n-1} \g{-w +\ds \f r n} = \lrr{2\pi}^{ \f{(n-1)}2}\ n^{ \f12 + nw}\  \lrc{\ds\f{\g{-nw}}{\g{-w}}}
\label{e11}
}
Observing that \h{\ds w+ 1- \f{n-r}{n} = w+\f r n} we have
\bea{
\ds\prod_{r=1}^{n-1} \Gamma\lrr{w+1-\f r n} = \ds\prod_{r=1}^{n-1} \Gamma\lrr{w+\f r n} =
(2\pi)^{\f{(n-1)}2} n^{\f12 - nw} \  \lrc{\ds\f{\g{nw}}{\g{w}}}
\label{e12}
}
Multiplying (\ref{e11}) and (\ref{e12}), and rearranging we get
\bea{
\ds\f{\Gamma(-nw)}{\Gamma(-w)} \aea \ds\lrc{\prod_{r=1}^{n-1}
\g{-w +\ds \f r n} \Gamma\lrr{w+1-\f r n}}\lrc{\f1{n\ (2\pi)^{n-1} }}\lrc{\ds\f{\Gamma(w)}{\Gamma(nw)}}
\label{eq13}
}
Observe that the right hand side of (\ref{eq13}) is well-defined for positive integer values of
\h{w} (since the singularities of the \h{\Gamma} function occur only at nonpositive integers).
Therefore
\bea{
\lim_{w\ra k} \ \ds\f{\Gamma(-nw)}{\Gamma(-w)} \aea \lrs{\lim_{w\ra k} \ds\lrc{\prod_{r=1}^{n-1}
\g{-w +\ds \f r n} \Gamma\lrr{w+1-\f r n}}}\lrc{\f1{n\ (2\pi)^{n-1} }}\lrc{\ds\f{\Gamma(k)}{\Gamma(nk)}}
\label{eq14}
}
Inserting (\ref{b1}) into (\ref{eq14}), and rewriting the limit in terms of \h{z=-w} we get
\bea{
\lim_{z\ra -k} \ \ds\f{\Gamma(nz)}{\Gamma(z)} \aea \ds\f{(-1)^k\ \Gamma(k)}{n^2 \ \Gamma(nk)}
}
as claimed \hfill \rule{7pt}{7pt}

\begin{figure}[h]
\centerline{\psfig{figure=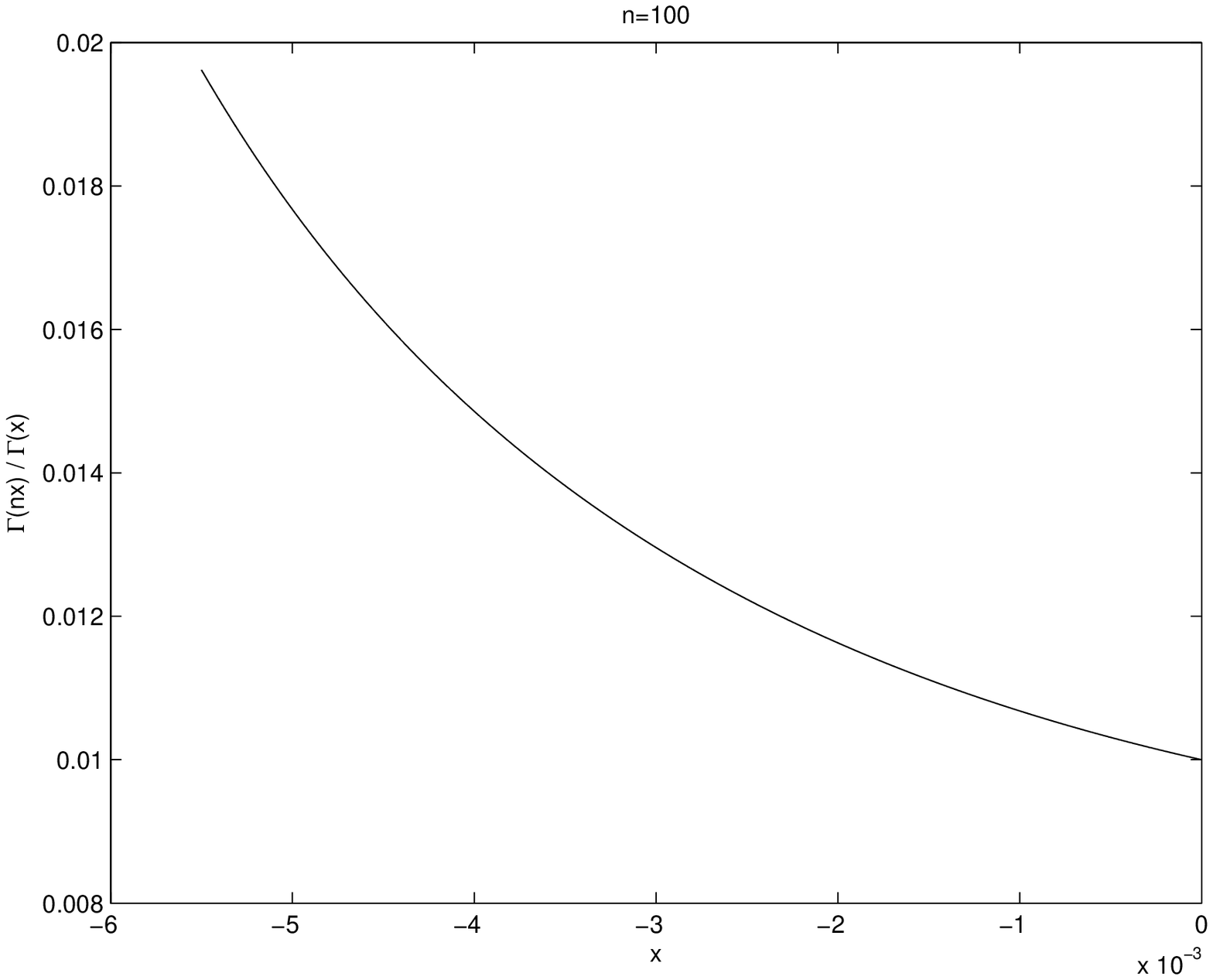,height=2in}\hspace{0.2in}\psfig{figure=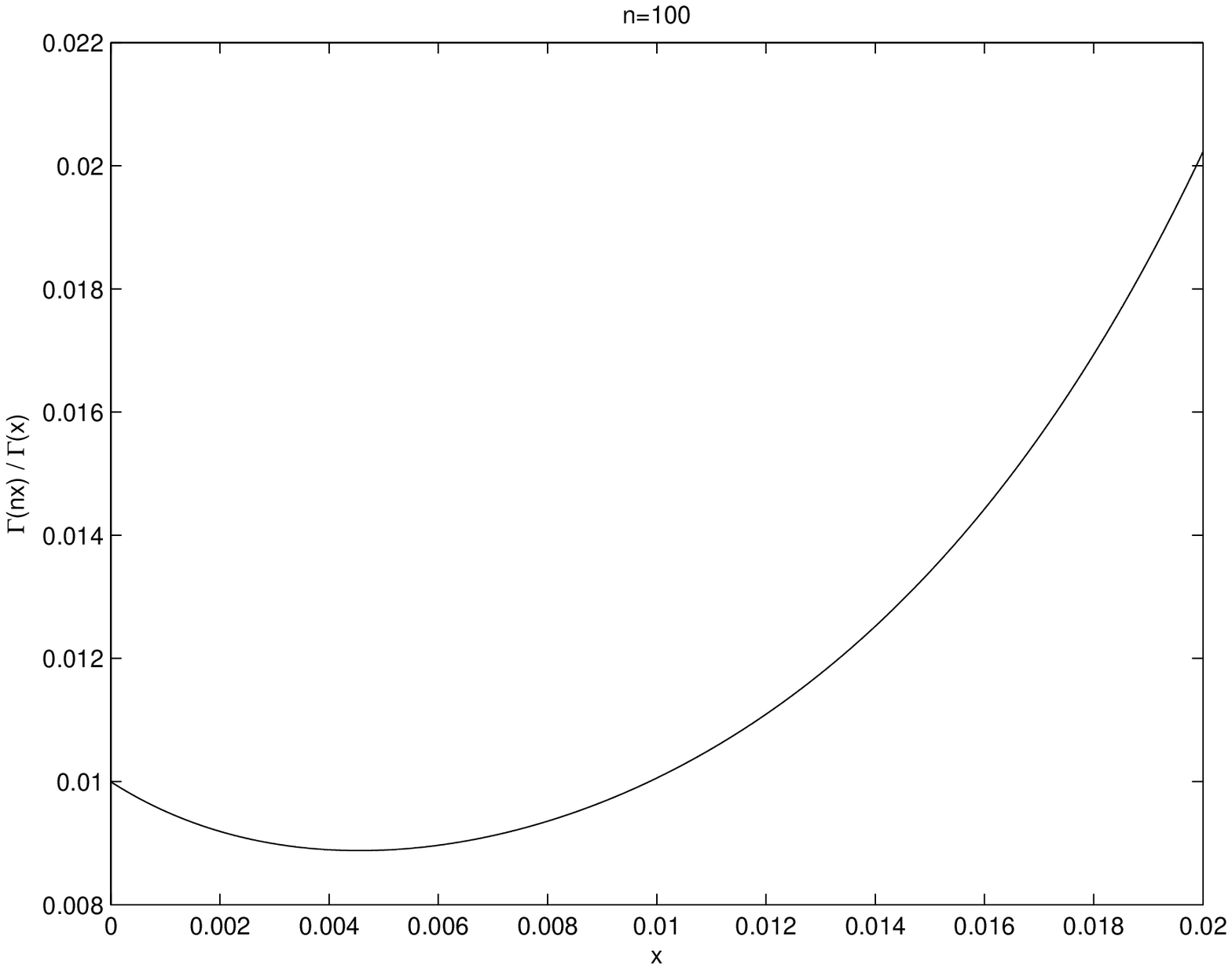,height=2in}}
\caption{Plot of \h{\ds\f{\Gamma(nz)}{\Gamma(z)}} to display the limits \h{z\ra 0^-} and \h{z\ra 0^+}}
\end{figure}

As an example, the function \h{\ds\f{\Gamma(nz)}{\Gamma(z)}} is plotted in Figure 1 to display the
limits \h{z\ra 0^-} and \h{z\ra 0^+}, for $n=100$.  The graphs converge to the predicted
value of 0.01.

\vspace{0.3in}

\appendix{\noindent \large\bf Appendix: Product of Sines}

\vspace{0.2in}

The following is an outline of the proof of (\ref{eq2}) as   presented in \cite{plan}.
For \h{n>1} and \h{0\leq k\leq n-1}, define
\h{\omega_k = e^{\f{2\pi i}{n} k}}. Expanding \h{z^n-1=
\prod_{k=1}^n (z-\omega_k)}, dividing by \h{z-1} and setting \h{z=1}, it follows that
\bea{
\prod_{i=1}^{n-1} \lrr{1-\omega_k} = n \label{a1}}
Using half-angle formula it follows that
\bea{
\left| 1-\omega_k\right| = 2\left|\sin\lrr{\f{\pi k}{n}}\right|\label{a2}
}
Setting \h{m=\lfloor n/2 \rfloor}, observing that
\h{\sin\lrr{\f{\pi k}n}> 0} for \h{1\leq k\leq n-1},
considering the odd and even \h{n}, and using (\ref{a1})
and (\ref{a2}) it follows that
\beas{
\ds\prod_{i=k}^m \sin^2\lrr{\f{\pi k}n} =
\prod_{k=1}^{n-1} \left|\sin\lrr{\f{\pi k}n}\right| =
\ds\f{n}{2^{n-1}}
}

\end{document}